\documentclass{amsart}
\usepackage{amsthm}
\usepackage{amsmath}
\usepackage{amsfonts}
\usepackage{amssymb}
\usepackage{graphicx}
\usepackage{float}
\usepackage{rotating}
\restylefloat{figure}
\usepackage{mathrsfs}
\usepackage{fancyhdr}
\usepackage[breaklinks]{hyperref}
\usepackage[numbers,sort&compress]{natbib}
\theoremstyle{plain}
\newtheorem{lem}{Lemma}

\newtheorem{thm}[lem]{Theorem}
\newtheorem{cor}{Corollary}

\theoremstyle{definition}

\newcommand{\R}{\mathbb R}

\newcommand{\N}{\mathbb N}
\newcommand{\dx}{\,dx}
\renewcommand{\d}{\,d}
\newcommand{\dw}{d}
\newcommand{\eps}{\varepsilon}
\renewcommand{\phi}{\varphi}
\newcommand{\norm}[1]{\Vert#1\Vert}

\renewcommand{\autoref}[1]{\text{Eq.}~\eqref{#1}}
\newcommand{\bea}{\begin{eqnarray}}
\newcommand{\eea}{\end{eqnarray}}
\newcommand{\beq}{\begin{equation}}
\newcommand{\eeq}{\end{equation}}
\begin{document}
\title{The two-component Camassa-Holm system in weighted $L_p$ spaces}
\author{Martin Kohlmann}
\address{Peter L. Reichertz Institute for Medical Informatics, University of Braunschweig, D-38106 Braunschweig, Germany}
\email{martin.kohlmann@plri.de}
\date{\today}
\begin{abstract} We present some new persistence results for the non-periodic two-component Camassa-Holm (2CH) system in weighted $L_p$ spaces. Working with moderate weight functions that are commonly used in time-frequency analysis, the paper generalizes some recent persistence results for the Camassa-Holm equation \cite{Br12} to its supersymmetric extension. As an application we discuss the spatial asymptotic profile of solutions to 2CH.
\end{abstract}
\maketitle
\tableofcontents
\section{Introduction and Main Results}\label{sec_intro}
As a model for the unidirectional irrotational free surface flow of a shallow layer of an inviscid fluid moving under the influence of gravity over a flat bed, the Camassa-Holm (CH) equation
\beq\label{CH} u_t - u_{txx} = uu_{xxx} + 2u_xu_{xx} - 3uu_x \eeq
has turned out to be suitable to capture typically nonlinear effects as wave breaking or peakons \cite{CH93,CM00}. The function $u(x,t)$ depends on a space variable $x\in\R$ and time $t\geq 0$ and represents the horizontal velocity of the fluid motion at a certain depth. Apart from its hydrodynamical relevance, the CH equation possesses various important mathematical properties: it has a bi-Hamiltonian structure \cite{CH93}, is obtained from Lagrange's variational principle \cite{IC07}, is completely integrable \cite{I05,I07} and allows for a geometric reformulation in terms of a geodesic flow on suitable diffeomorphism groups \cite{L07}.

Several two-component generalizations of the CH equation have been introduced recently and the most popular of them reads
\beq\label{2CHsmooth}
\left\{
\begin{array}{rl}
u_t   - u_{txx}    & = uu_{xxx} + 2u_xu_{xx} - 3uu_x -\rho\rho_x, \\
\rho_t & =-(\rho u)_x.
\end{array}
\right.
\eeq
This system first appeared in \cite{OR96} and had been discussed in many further publications: The associated Cauchy problem in the real line setting has been subject of, e.g., \cite{ELY07,GL11,GY10,GZ10}. Conservative solutions on the full line are presented in \cite{GHR12} and global weak solutions have been studied in \cite{GY11,GY12}. The inverse scattering method has been applied in \cite{HI11} to \eqref{2CHsmooth} where the authors also show that 2CH can be obtained from the following Lax pair \cite{CI08,LL09}
\begin{align}
\psi_{xx} & = \left(-\lambda^2\rho^2+\lambda m+\frac{1}{4}\right)\psi, \nonumber\\
\psi_t    & = \left(\frac{1}{2\lambda}-u\right)\psi_x + \frac{1}{2}u_x\psi; \nonumber
\end{align}
here, $m=u-u_{xx}$. Furthermore, 2CH is bi-Hamiltonian with the Hamiltonians
$$H_1=\frac{1}{2}\int (um+\rho^2)\dx,\qquad H_2=\frac{1}{2}\int(u\rho^2+u^3+uu_x^2)\dx,$$
cf., e.g., \cite{CI08,LL09,P06}, and of variational nature for the Lagrangian $H_1$. Concerning geometric aspects of \eqref{2CHsmooth} it is shown in \cite{GO06} that 2CH reexpresses geodesic motion on certain semidirect product Lie groups. It is explained in \cite{CI08} that the short-wave limit of the 2CH system (leading to the two-component Hunter-Saxton equation, in fact) possesses the peakon solutions
\begin{align}
m(x,t) & = \sum_{k=1}^N m_k(t)\delta(x-x_k(t)), \nonumber\\
u(x,t) & = -\frac{1}{2}\sum_{k=1}^N m_k(t)|x-x_k(t)|, \nonumber\\
\rho(x,t) & = \sum_{k=1}^N\rho_k(t)\theta(x-x_k(t)).\nonumber
\end{align}
Some explicit examples like traveling waves and multi-kink solutions of \eqref{2CHsmooth} are discussed in \cite{CLZ06,LL09,P06}. That 2CH can be seen as a supersymmetric generalization of the CH equation is explained in \cite{LL09,P06}. As it is of particular importance for the issues of the present paper, we shall comment on persistence properties and the asymptotic behavior of solutions to the 2CH system: In \cite{H09}, the author proves that solutions of 2CH have infinite propagation speed in the sense that (non-zero) compactly supported initial data lead to solutions that are not compactly supported for any positive time. Instead, these solutions possess exponentially decaying profiles for large values of the spatial variable. The results of \cite{H09} generalize the work presented in \cite{C05,H05,H08,H09b} on infinite propagation speed and asymptotic profiles for the CH equation. More details on persistence and unique continuation and asymptotic profiles for the 2CH system can be found in \cite{G12,GN11}.

It is an immediate consequence of Kato's semigroup theory \cite{K75} that \eqref{2CHsmooth} is well-posed in the spaces $H^s\times H^{s-1}$ for $s>5/2$, where $H^s$ is the Sobolev space of order $s\geq 0$ on the real line. Precisely, given $z_0=(u_0,\rho_0)\in H^s\times H^{s-1}$ for some $s>5/2$, there exists a maximal number $T>0$ and a unique solution $z=(u,\rho)\in C([0,T);H^s\times H^{s-1})\cap C^1([0,T);H^{s-1}\times H^{s-2})$ to \eqref{2CHsmooth} with $z(0)=z_0$ such that the mapping $z_0\mapsto z$ is continuous. In this paper, we intend to find a large class of weight functions $\phi$ such that
$$\sup_{t\in[0,T)}\left(\norm{\phi u(t)}_p+\norm{\phi u_x(t)}_p+\norm{\phi u_{xx}(t)}_p+\norm{\phi\rho(t)}_p+\norm{\phi \rho_x(t)}_p\right)<\infty$$
where $\norm{\cdot}_p$ denotes the usual $L_p$ norm. This way we obtain a persistence result on solutions $z$ to \eqref{2CHsmooth} in the weighted $L_p$ spaces $L_{p,\phi}:=L_p(\R,\phi^pdx)$. As a consequence and an application we determine the spatial asymptotic behavior of certain solutions to 2CH. Our results generalize the work of \cite{Br12} on persistence and non-persistence of solutions to \eqref{CH} in $L_{p,\phi}$.\ We will work with moderate weight functions which appear with regularity in the theory of time-frequency analysis \cite{AG01,G07} and have led to optimal results for the CH equation in \cite{Br12}. In Section~\ref{sec_prelim} we will characterize a subclass of moderate weight functions to be admissible for 2CH. Then our main theorem reads as follows:
\thm\label{thm1} Let $s>3$ and $2\leq p\leq\infty$. Let $z=(u,\rho)\in C([0,T), H^s\times H^{s-1})\cap C^1([0,T), H^{s-1}\times H^{s-2})$ be the strong solution to \eqref{2CHsmooth} starting from $z_0=(u_0,\rho_0)$ so that $\phi u_0,\phi u_{0x},\phi u_{0xx},\phi\rho_0,\phi\rho_{0x}\in L_p(\R)$ for an admissible weight function $\phi$ of 2CH. Let
$$M:=\sup_{t\in[0,T)}\left\{\norm{u(t)}_\infty+\norm{u_x(t)}_\infty+\norm{u_{xx}(t)}_\infty+\norm{\rho(t)}_\infty+\norm{\rho_x(t)}_\infty\right\}.$$
Then there is a constant $C>0$ depending only on the weight $\phi$ such that
\begin{align}&\norm{\phi u(t)}_p+\norm{\phi u_x(t)}_p+\norm{\phi u_{xx}(t)}_p+\norm{\phi\rho(t)}_p+\norm{\phi \rho_x(t)}_p\nonumber\\
&\qquad\leq e^{CMt}\left(\norm{\phi u_0}_p+\norm{\phi u_{0x}}_p+\norm{\phi u_{0xx}}_p+\norm{\phi\rho_0}_p+\norm{\phi\rho_{0x}}_p\right),\nonumber\end{align}
for all $t\in[0,T)$.
\endthm\rm
The basic class of examples for admissible weight functions is given by the family of functions
$$\phi_{a,b,c,d}(x) = e^{a|x|^b}(1+|x|)^c\left(\log(e+|x|)\right)^d,$$
where we demand $a\geq 0$, $0\leq b\leq 1$ and $ab<1$, cf.~Section \ref{sec_prelim}. Let us comment on two particular results:
\begin{itemize}
\item Let $\phi=\phi_{0,0,c,0}$, $c>0$, and $p=\infty$ in Theorem~\ref{thm1}. For this choice, Theorem~\ref{thm1} says that the algebraic decay of the initial datum $z_0=(u_0,\rho_0)$,
$$|u_0(x)|+|u_{0x}(x)|+|u_{0xx}(x)|+|\rho_0(x)|+|\rho_{0x}(x)|\leq C(1+|x|)^{-c},$$
for all $x\in\R$, is preserved by the solution $z=(u,\rho)$ with $z(0)=z_0$ on $[0,T)$, i.e.,
$$
|u(x,t)|+|u_{x}(x,t)|+|u_{xx}(x,t)|+|\rho(x,t)|+|\rho_{x}(x,t)|\leq C'(1+|x|)^{-c},
$$
for all $(x,t)\in\R\times[0,T)$, where $C,C'>0$ are constants. Theorem~\ref{thm1} thus generalizes the main result of Ni and Zhou \cite{NZ12} on algebraic decay rates of strong solutions to the CH equation.
\item Let $\phi=\phi_{a,1,0,0}1_{\{x\geq 0\}}+1_{\{x<0\}}$, $0\leq a<1$. By our definitions in Section~\ref{sec_prelim} $\phi$ is an admissible weight function for 2CH. Let furthermore $p=\infty$ in Theorem~\ref{thm1}. Then one deduces that 2CH preserves the pointwise decay $O(e^{-ax})$ of its solutions as $x\to\infty$, for any $t>0$. Analogously one concludes that, for $x\to-\infty$, the decay $O(e^{ax})$ is preserved during the evolution. A corresponding result on persistence of strong solutions of the CH equation can be found in Theorem 1.2 of \cite{HMPZ07}; it is worthwhile to note that Himonas, Misio{\l}ek, Ponce and Zhou also use weight functions and Gronwall's Lemma to obtain a proof of \cite[Theorem 1.2]{HMPZ07}.
\end{itemize}
Theorem~\ref{thm1} does not apply to the family of weights $\phi_{1,1,c,d}$. In the following corollary however we may choose $\phi=\phi_{1,1,c,d}$ if $c<0$, $d\in\R$ and $\frac{1}{|c|}<p\leq\infty$. For the notion of $v$-moderate weights we refer the reader to Section~\ref{sec_prelim}.
\cor\label{cor1} Let $p\leq 2\leq\infty$ and let $\phi\colon\R\to(0,\infty)$ be a locally absolutely continuous and $v$-moderate weight function satisfying $|\phi'(x)|\leq A|\phi(x)|$ a.e., for some $A>0$, $\inf v>0$ and $ve^{-|\cdot|}\in L_p(\R)$. Assume that $\phi u_0$, $\phi u_{0x}$, $\phi u_{0xx}$, $\phi\rho_0$, $\phi\rho_{0x}\in L_p(\R)$ and $\phi^{1/2} u_0$, $\phi^{1/2} u_{0x}$, $\phi^{1/2} u_{0xx}$, $\phi^{1/2}\rho_0$, $\phi^{1/2}\rho_{0x}\in L_2(\R)$. For $s>3$, let $z=(u,\rho)\in C([0,T), H^s\times H^{s-1})\cap C^1([0,T), H^{s-1}\times H^{s-2})$ be the strong solution to \eqref{2CHsmooth} starting from $z_0=(u_0,\rho_0)$. Then
$$\sup_{t\in[0,T)}\left(\norm{\phi u(t)}_p + \norm{\phi u_x(t)}_p + \norm{\phi u_{xx}(t)}_p + \norm{\phi\rho(t)}_p + \norm{\phi\rho_x(t)}_p\right)$$
and
$$\sup_{t\in[0,T)}\bigg(\norm{\phi^{1/2}u(t)}_2 + \norm{\phi^{1/2}u_x(t)}_2 + \norm{\phi^{1/2}u_{xx}(t)}_2 \nonumber\\
+ \norm{\phi^{1/2}\rho(t)}_2 + \norm{\phi^{1/2}\rho_x(t)}_2\bigg)$$
are finite.
\endcor\rm
For the particular choice $c=d=0$ and $p=\infty$, we conclude from
$$|u_0(x)| + |u_{0x}(x)| + |u_{0xx}(x)| + |\rho_0(x)| + |\rho_{0x}(x)| \leq Ce^{-|x|}$$
for any $x\in\R$ that the unique solution $z=(u,\rho)\in C([0,T);\,H^s\times H^{s-1})$ of 2CH with $z(0)=(u_0,\rho_0)$ satisfies
$$|u(x,t)| + |u_x(x,t)| + |u_{xx}(x,t)| + |\rho(x,t)| + |\rho_x(x,t)| \leq C'e^{-|x|}$$
on $\R\times[0,T)$. In the following corollary we compute the spatial asymptotic profiles of solutions with exponential decay. As a further consequence we may infer that the peakon-like decay $O(e^{-|x|})$ mentioned above is the fastest possible decay for a nontrivial solution $z$ of 2CH to propagate.
\cor\label{cor2} Let $\psi(x)=e^{|x|/2}(1+|x|)^{1/2}\left(\log(e+|x|)\right)^d$, for some $d>1/2$. For $s>3$, let $z_0=(u_0,\rho_0)\in H^s\times H^{s-1}$ be nonzero and assume that
\beq\sup_{x\in\R}\left\{\psi(x)\left(|u_0(x)|+|u_{0x}(x)|+|u_{0xx}(x)|+|\rho_0(x)|+|\rho_{0x}(x)|\right)\right\}<\infty.\label{cond}\eeq
Then condition \eqref{cond} is conserved for the solution $z=(u,\rho)\in C([0,T), H^s\times H^{s-1})\cap C^1([0,T), H^{s-1}\times H^{s-2})$ of 2CH starting from $z_0$, and we have the asymptotic behavior
$$z(t)=
\left\{
\begin{array}{lll}
  z_0 + e^{-x}t\begin{pmatrix}\Phi^+(t)+\eps^+_1(x,t) \\ \eps^+_2(x,t)\end{pmatrix},&\text{as}& x\to\infty, \\
  z_0 - e^{x}t\begin{pmatrix}\Phi^-(t)+\eps^-_1(x,t) \\ \eps^-_2(x,t)\end{pmatrix},&\text{as}& x\to-\infty,
\end{array}
\right.
$$
for all $t\in[0,T)$, where $\eps^+_{1,2}(x)\to 0$ as $x\to\infty$, $\eps^-_{1,2}(x)\to 0$ as $x\to-\infty$ and $c_1\leq\Phi^{\pm}(t)\leq c_2$ with $c_1,c_2>0$ independent of $t$.
\endcor\rm
Let us compare our results with the basic paper \cite{HMPZ07} on persistence properties and unique continuation of solutions for the Camassa-Holm equation. Theorem~\ref{thm1} on persistence of solutions of the 2CH system generalizes Theorem 1.2 in \cite{HMPZ07} dealing with persistence of solutions for the CH equation. Corollary~\ref{cor2} about the asymptotic profile for solutions of 2CH corresponds to Theorem 1.4 in \cite{HMPZ07} about the asymptotic profile of solutions to the CH equation. Corollary~\ref{cor2} also corresponds to Theorem 3.5 of \cite{H09} dealing with infinite propagation speed for 2CH and the spatial asymptotic profile obtained for compactly supported initial data.

Our paper is organized as follows: In Section~\ref{sec_prelim} we present some fundamentals concerning moderate weight functions and the functional analytic setting for 2CH. In Section~\ref{sec_proofs} we prove Theorem~\ref{thm1} and its two corollaries.\\[.25cm]
\textbf{Acknowledgement.} The author is grateful to the referees whose suggestions helped to improve the first version of the paper.
\section{Preliminaries}\label{sec_prelim}
The formulation \eqref{2CHsmooth} of the 2CH only holds for smooth solutions. Therefore we first rewrite \eqref{2CHsmooth} as
\beq
\left\{
\begin{array}{ll}
u_t   +uu_x    & = P(D)(u^2+\frac{1}{2}u_x^2+\frac{1}{2}\rho^2), \\
\rho_t+u\rho_x & =-\rho u_x,
\end{array}
\right.\quad P(D)=-\partial_x(1-\partial_x^2)^{-1}.
\label{2CH}
\eeq
The Sobolev spaces $H^s=H^s(\R)$, $s\geq 0$, on the line are the Hilbert spaces
$$H^s = \left\{u\in L_2(\R);\,\int_\R|\hat u(x)|^2(1+x^2)^s\dx<\infty\right\}$$
equipped with the norms
$$\norm{u}_{H^s}^2=\int_\R|\hat u(x)|^2(1+x^2)^s\dx.$$
The pseudodifferential operator $(1-\partial_x^2)^{-1}$ has the symbol $\frac{1}{1+x^2}$ and hence defines an isomorphism $H^s\to H^{s+2}$ for any $s\geq 0$. Moreover, $(1-\partial_x^2)^{-1}f(x)=(G*f)(x)$ with the kernel $G(x)=\frac{1}{2}e^{-|x|}$. Some standard computations show that
\beq\label{Gx}(\partial_x G)(x) = -\frac{1}{2}e^{-x}1_{\{x\geq 0\}} + \frac{1}{2}e^x1_{\{x<0\}} = -\frac{1}{2}\text{sign}(x)e^{-|x|}\eeq
in the weak sense and that
\beq\label{Gxx}\partial_x^2G=G-\delta\eeq
with the Dirac distribution $\delta$. It follows that $P(D)$ is a $\psi$do with the symbol $-\frac{\text{i}x}{1+x^2}$ and the Green's function $-\partial_x G$, and that $P(D)\colon H^s\to H^{s+1}$ is an isomorphism for any $s\geq 0$. By Sobolev's embedding theorem, we also have $H^s(\R)\subset C^k_b(\R)$, $s>k+1/2$, where $C^k_b(\R)$ is the space of $k$ times continuously differentiable functions on $\R$ which are, together with any of their derivatives, bounded with respect to $\norm{\cdot}_\infty$. In particular, for $(u,v)\in H^{s}\times H^{s-1}$, $s>5/2$, the functions $u,u_x,u_{xx},v,v_x$ are continuous and bounded on $\R$.

A function $v\colon\R^n\to\R$ is called sub-multiplicative if $v(x+y)\leq v(x)v(y)$, for all $x,y\in\R^n$. Let $v$ be a sub-multiplicative function. A positive function $\phi$ on $\R^n$ is called $v$-moderate if there exists a constant $c>0$ such that $\phi(x+y)\leq cv(x)\phi(y)$, for all $x,y\in\R^n$. We say that $\phi$ is moderate if it is $v$-moderate for some sub-multiplicative function $v$. As shown in \cite{Br12}, a positive function $\phi$ is $v$-moderate with constant $c$ if and only if for any two measurable functions $f_1$ and $f_2$ and $1\leq p\leq \infty$ the weighted Young inequality
\beq\label{young}\norm{(f_1*f_2)\phi}_p\leq c\norm{f_1v}_1\norm{f_2\phi}_p\eeq
holds.

We say that $\phi\colon\R\to(0,\infty)$ is an admissible weight function for the 2CH equation if it is locally absolutely continuous such that for some $A>0$ we have a.e.\ $|\phi'(x)|\leq A|\phi(x)|$ and $\phi$ is $v$-moderate with a sub-multiplicative function $v$ satisfying $\inf v>0$ and
\beq\label{admissible}\int_\R v(x)e^{-|x|}\dx<\infty.\eeq
It is easily checked that the functions $\phi_{a,b,c,d}$ presented in Section~\ref{sec_intro} are admissible for 2CH if $a\geq 0$, $0\leq b\leq 1$ and $ab<1$; see \cite{Br12} for further details.
\section{Proofs}\label{sec_proofs}
\subsection{Proof of Theorem~\ref{thm1}}
Let $F(u,\rho)=u^2+\frac{1}{2}u_x^2+\frac{1}{2}\rho^2$ and assume that $\phi$ is $v$-moderate satisfying the conditions specified in Section~\ref{sec_prelim}. Our first observation is that the first row equation in \eqref{2CH} can be rewritten as $u_t+uu_x+\partial_xG*F(u,\rho)=0$. For any $n\in\N$, let $\phi_n(x)=\min\{\phi(x),n\}$. Then $\phi_n\colon\R\to\R$ is locally absolutely continuous, $\norm{\phi_n}_\infty\leq n$ and $|\phi_n'(x)|\leq A|\phi_n(x)|$ a.e.\ on $\R$.
Moreover, as shown in \cite{Br12}, the $n$-truncations $\phi_n$ are again $v$-moderate. Let $p\in[2,\infty)$. We multiply the first row equation of \eqref{2CH} by $\phi_n|\phi_nu|^{p-2}\phi_nu$ and integrate to obtain
\begin{align}&\int_\R(\phi_nu)_t\phi_nu|\phi_nu|^{p-2}\dx + \int_\R |\phi_nu|^pu_x\dx \nonumber\\ &\qquad + \int_\R \phi_n(\partial_xG*F(u,\rho)) |\phi_nu|^{p-2}\phi_nu\dx=0.\label{inteq1}\end{align}
We denote the three terms on the left-hand side of \eqref{inteq1} as $\mathcal I_1$, $\mathcal I_2$ and $\mathcal I_3$ and observe that
$$
\mathcal I_1 = \frac{1}{p}\int_\R\frac{\dw}{\dw t}\left[\sqrt{(\phi_nu)^2}\right]^p\dw x = \frac{1}{p}\frac{\dw}{\dw t}\norm{\phi_nu}_p^p
=\norm{\phi_nu}_p^{p-1}\frac{d}{dt}\norm{\phi_nu}_p,
$$
that $|\mathcal I_2|\leq M\norm{\phi_nu}_p^p$ and that, by H\"older's inequality,
$$|\mathcal I_3| \leq \norm{\phi_n(\partial_xG*F(u,\rho))}_p\norm{(\phi_nu)^{p-1}}_{\frac{1}{1-1/p}}=
\norm{\phi_n(\partial_xG*F(u,\rho))}_p\norm{\phi_nu}_{p}^{p-1}.$$
Hence \eqref{inteq1} yields
\beq\label{aux1}\frac{d}{dt}\norm{\phi_nu}_p\leq M\norm{\phi_nu}_p + \norm{\phi_n(\partial_xG*F(u,\rho))}_p,\eeq
and using \eqref{Gx} and \eqref{young} we find
\begin{align}
\norm{\phi_n(\partial_xG*F(u,\rho))}_p
& \leq C_1\norm{(\partial_xG)v}_1\norm{\phi_nF(u,\rho)}_p\nonumber\\
& \leq C_2\norm{\phi_nF(u,\rho)}_p,\label{estimateGx}
\end{align}
where $C_1$ and $C_2$ are independent of $n$. As
\begin{align}
\norm{\phi_nF(u,\rho)}_p & \leq\norm{\phi_nu^2}_p+\frac{1}{2}\norm{\phi_nu_x^2}_p+\frac{1}{2}\norm{\phi_n\rho^2}_p \nonumber\\
& \leq M\left(\norm{\phi_nu}_p + \norm{\phi_nu_x}_p + \norm{\phi_n\rho}_p\right)\label{phinF}
\end{align}
we conclude from \eqref{aux1}, \eqref{estimateGx} and \eqref{phinF} that
\beq\label{u}\frac{d}{dt}\norm{\phi_nu}_p\leq (C_2+1)M(\norm{\phi_nu}_p + \norm{\phi_nu_x}_p + \norm{\phi_n\rho}_p).\eeq
Multiplying the identity $u_{tx}+uu_{xx}+u_x^2+\partial_x^2G*F(u,\rho)=0$ with $\phi_n|\phi_nu_x|^{p-2}\phi_nu_x$ we obtain as before
\begin{align}\label{inteq1x}
&\frac{1}{p}\frac{d}{dt}\norm{\phi_nu_x}_p^p + \int_\R uu_{xx}\phi_n|\phi_nu_x|^{p-2}\phi_nu_x\dx+\int_\R u_{x}^2\phi_n|\phi_nu_x|^{p-2}\phi_nu_x\dx
\nonumber\\ &\qquad +\int_\R\phi_n(\partial_x^2G*F(u,\rho))|\phi_nu_x|^{p-2}\phi_nu_x\dx = 0.
\end{align}
Let $\mathcal I_4$, $\mathcal I_5$ and $\mathcal I_6$ denote the three integrals on the left-hand side of (\ref{inteq1x}). As before, we will work with the estimates $|\mathcal I_4|\leq M\norm{\phi_nu}_p\norm{\phi_nu_x}_p^{p-1}$, $|\mathcal I_5|\leq M\norm{\phi_nu_x}_p^p$ and $|\mathcal I_6|\leq\norm{\phi_n(\partial_x^2G*F(u,\rho))}_p\norm{\phi_nu_x}_p^{p-1}$. Thus
\beq\label{aux2}\frac{d}{dt}\norm{\phi_nu_x}_p \leq M\left(\norm{\phi_nu}_p+\norm{\phi_nu_x}_p\right) + \norm{\phi_n(\partial_x^2G*F(u,\rho))}_p.\eeq
Using once again that $\norm{vG}_1<\infty$, \eqref{Gxx} and (\ref{phinF}) we get
$$\norm{\phi_n(\partial_x^2G*F(u,\rho))}_p\leq C_3\norm{\phi_nF(u,\rho)}_p\leq C_3M(\norm{\phi_nu}_p + \norm{\phi_nu_x}_p + \norm{\phi_n\rho}_p).$$
This achieves
\beq\label{ux}\frac{d}{dt}\norm{\phi_nu_x}_p\leq (C_3+1)M(\norm{\phi_nu}_p + \norm{\phi_nu_x}_p + \norm{\phi_n\rho}_p).\eeq
Multiplying the equation $u_{txx}+uu_{xxx}+3u_xu_{xx}+\partial_xG*F(u,\rho)-\partial_xF(u,\rho)=0$ with $\phi_n|\phi_nu_{xx}|^{p-2}\phi_nu_{xx}$ yields
\begin{align}\label{inteq1xx}
&\frac{1}{p}\frac{d}{dt}\norm{\phi_nu_{xx}}_p^p + \int_\R uu_{xx}u_{xxx}\phi_n^2|\phi_nu_{xx}|^{p-2}\dx \nonumber\\
& \quad + 3\int_\R u_xu_{xx}^2\phi_n^2|\phi_nu_{xx}|^{p-2}\dx + \int_\R\phi_n(\partial_x G*F(u,\rho))|\phi_nu_{xx}|^{p-2}\phi_n u_{xx}\dx \nonumber\\
& \quad - \int_\R\phi_n\partial_xF(u,\rho)|\phi_nu_{xx}|^{p-2}\phi_nu_{xx}\dx = 0.
\end{align}
To estimate the four integrals $\mathcal I_7,\ldots,\mathcal I_{10}$ on the left-hand side of \eqref{inteq1xx}, we note first that $|\mathcal I_8|\leq 3 M\norm{\phi_nu_{xx}}_p^p$ and that $|\mathcal I_9|\leq\norm{\phi_n(\partial_x G*F(u,\rho))}_p\norm{\phi_nu_{xx}}_p^{p-1}$ and $|\mathcal I_{10}|\leq\norm{\phi_n\partial_xF(u,\rho)}_p\norm{\phi_nu_{xx}}_p^{p-1}$. Similarly to \eqref{phinF}, we estimate
\begin{align}\label{phinFx}\norm{\phi_n\partial_xF(u,\rho)}_p&\leq 2\norm{\phi_nuu_x}_p+\norm{\phi_nu_xu_{xx}}_p+\norm{\phi_n\rho\rho_{x}}_p\nonumber\\
&\leq 2M\left(\norm{\phi_nu}_p + \norm{\phi_nu_x}_p + \norm{\phi_n\rho}_p\right),\end{align}
so that, by \eqref{estimateGx}, \eqref{phinF} and \eqref{phinFx},
$$|\mathcal I_9|,|\mathcal I_{10}|\leq C_4M\left(\norm{\phi_nu}_p + \norm{\phi_nu_x}_p + \norm{\phi_n\rho}_p\right)\norm{\phi_nu_{xx}}_p^{p-1}.$$
The integral $\mathcal I_7$ may be rewritten as
$$\mathcal I_7 = \int_\R\left[\partial_x(\phi_nu_{xx})-\phi_{nx}u_{xx}\right]uu_{xx}\phi_n|\phi_nu_{xx}|^{p-2}\dx=\mathcal I_7^{(1)}-\mathcal I_7^{(2)}.$$
As $|\phi_{nx}|\leq A|\phi_n|$, we have $|\mathcal I_7^{(2)}|\leq AM\norm{\phi_nu_{xx}}_p^p$. For any function $g\colon\R\times[0,T)\to\R$ which is weakly differentiable in the first argument one has
\beq\label{auxint}\frac{\partial}{\partial x}\left(\frac{|g(x,t)|^p}{p}\right)=|g(x,t)|^{p-2}g(x,t)g_x(x,t).\eeq
Since $s>3$ and $\phi_n$ is locally absolutely continuous, \autoref{auxint} holds a.e.\ for $g=\phi_nu_{xx}$ and we have
$$\mathcal I_7^{(1)}=\int_\R u\partial_x\left(\frac{|\phi_nu_{xx}|^p}{p}\right)\dx=-\frac{1}{p}\int_\R u_x|\phi_nu_{xx}|^p\dx;$$
the boundary terms vanish when performing integration by parts in view of Sobolev's embedding theorem. Hence $|\mathcal I_7^{(1)}|\leq\frac{M}{p}\norm{\phi_nu_{xx}}_p^p$. We have shown that
\begin{align}
\label{aux3}\frac{d}{dt}\norm{\phi_nu_{xx}}_p &\leq (A+4)M\norm{\phi_nu_{xx}}_p + 2M\left(\norm{\phi_nu}_p + \norm{\phi_nu_x}_p + \norm{\phi_n\rho}_p\right)\nonumber\\
&\qquad+ \norm{\phi_n(\partial_xG*F(u,\rho))}_p \\
\label{uxx}&\leq C_5M\left(\norm{\phi_nu}_p+\norm{\phi_nu_x}_p+\norm{\phi_nu_{xx}}_p+\norm{\phi_n\rho}_p\right).
\end{align}
We now multiply the equation $\rho_t+u\rho_x+\rho u_x=0$ with $\phi_n|\phi_n\rho|^{p-2}\phi_n\rho$ and integrate to obtain the identity
\beq\label{inteq2}\frac{1}{p}\frac{d}{dt}\norm{\phi_n\rho}_p^p+\int_\R\rho\rho_xu\phi_n^2|\phi_n\rho|^{p-2}\dx+\int_\R(\phi_n\rho)^2u_x|\phi_n\rho|^{p-2}\dx=0.\eeq
Let $\mathcal J_1$ and $\mathcal J_2$ be the integrals on the left-hand side of \eqref{inteq2}. As above, we infer that $|\mathcal J_1|\leq M\norm{\phi_n\rho_x}_p\norm{\phi_n\rho}_p^{p-1}$ and that $|\mathcal J_2|\leq M\norm{\phi_n\rho}_p^p$. This yields
\beq\label{rho}\frac{d}{dt}\norm{\phi_n\rho}_p\leq M(\norm{\phi_n\rho}_p+\norm{\phi_n\rho_x}_p).\eeq
Multiplying the identity $\rho_{tx}+\rho_{xx}u+2\rho_xu_x+\rho u_{xx}=0$ with $\phi_n|\phi_n\rho_x|^{p-2}\phi_n\rho_x$ and integrating we also obtain
\begin{align}
&\frac{1}{p}\frac{d}{dt}\norm{\phi_n\rho_x}_p^p+\int_\R\rho_{xx}\rho_xu\phi_n^2|\phi_n\rho_x|^{p-2}\dx+2\int_\R\rho_x^2u_x\phi_n^2|\phi_n\rho_x|^{p-2}\dx\nonumber\\
&\qquad+\int_\R\rho\rho_xu_{xx}\phi_n^2|\phi_n\rho_x|^{p-2}\dx=0.\label{inteq2x}
\end{align}
For the three integrals $\mathcal J_3$, $\mathcal J_4$ and $\mathcal J_5$ on the left-hand side of \eqref{inteq2x} we have $|\mathcal J_5|\leq M\norm{\phi_n\rho}_p\norm{\phi_n\rho_x}_p^{p-1}$, $|\mathcal J_4|\leq 2M\norm{\phi_n\rho_x}_p^p$ and $\mathcal J_3$ may be decomposed as
$$
\mathcal J_3 = \int_\R\left[\partial_x(\phi_n\rho_x)-\phi_{nx}\rho_x\right]u|\phi_n\rho_x|^{p-2}\phi_n\rho_x\dx=\mathcal J_3^{(1)}-\mathcal J_3^{(2)}.
$$
Applying \eqref{auxint} to $g=\phi_n\rho_x$ we conclude
$$\mathcal J_3^{(1)}=\int_\R u\partial_x\left(\frac{|\phi_n\rho_x|^p}{p}\right)\dx=-\frac{1}{p}\int_\R u_x|\phi_n\rho_x|^p\dx;$$
again the boundary terms vanish when performing integration by parts in view of Sobolev's embedding theorem. Hence $|\mathcal J_3^{(1)}|\leq\frac{M}{p}\norm{\phi_n\rho_x}_p^p$. Using that $|{\phi_{nx}}|\leq A|\phi_n|$ it is clear that $|\mathcal J_3^{(2)}|\leq AM\norm{\phi_n\rho_x}_p^p$. Thus
\beq\label{rhox}\frac{d}{dt}\norm{\phi_n\rho_x}_p\leq M\norm{\phi_n\rho}_p + (3+A)M\norm{\phi_n\rho_x}_p.\eeq
Combining Eqs.\ \eqref{u}, \eqref{ux}, \eqref{uxx}, \eqref{rho} and \eqref{rhox}, there exists a constant $C$, only depending on $A$ and $v$, such that
\begin{align}
& \frac{d}{dt}\left(\norm{\phi_nu}_p+\norm{\phi_nu_x}_p+\norm{\phi_nu_{xx}}_p+\norm{\phi_n\rho}_p+\norm{\phi_n\rho_x}_p\right) \nonumber\\
& \qquad\leq CM\left(\norm{\phi_nu}_p+\norm{\phi_nu_x}_p+\norm{\phi_nu_{xx}}_p+\norm{\phi_n\rho}_p+\norm{\phi_n\rho_x}_p\right), \nonumber
\end{align}
so that, by Gronwall's Lemma,
\begin{align}
& \norm{\phi_nu(t)}_p+\norm{\phi_nu_x(t)}_p+\norm{\phi_nu_{xx}(t)}_p+\norm{\phi_n\rho(t)}_p+\norm{\phi_n\rho_x(t)}_p \nonumber\\
& \quad\leq \left(\norm{\phi_nu_0}_p+\norm{\phi_nu_{0x}}_p+\norm{\phi_nu_{0xx}}_p+\norm{\phi_n\rho_0}_p+\norm{\phi_n\rho_{0x}}_p\right)e^{CMt},\forall t\in[0,T).\nonumber
\end{align}
Since $\phi_n(x)\to\phi(x)$ a.e.\ as $n\to\infty$ and $\phi u_0,\phi u_{0x},\phi u_{0xx},\phi\rho_0,\phi\rho_{0x}\in L_p(\R)$ the assertion of the theorem follows for the case $p\in[2,\infty)$. Since $\norm{\cdot}_\infty=\lim_{p\to\infty}\norm{\cdot}_p$ it is clear that the theorem also applies for $p=\infty$.\hfill$\square$
\subsection{Proof of Corollary~\ref{cor1}}
As explained in \cite{Br12}, the function $\phi^{1/2}$ is a $v^{1/2}$-moderate weight satisfying $|(\phi^{1/2})'(x)|\leq\frac{A}{2}\phi^{1/2}(x)$, $\inf v^{1/2}>0$ and $v^{1/2}e^{-|\cdot|}\in L_1(\R)$. We apply Theorem~\ref{thm1} with $p=2$ to the weight $\phi^{1/2}$ and obtain
\begin{align}\label{est.5}
&\norm{\phi^{1/2}u(t)}_2 + \norm{\phi^{1/2}u_x(t)}_2 + \norm{\phi^{1/2}u_{xx}(t)}_2 + \norm{\phi^{1/2}\rho(t)}_2 + \norm{\phi^{1/2}\rho_x(t)}_2\\
&\quad\leq\left(\norm{\phi^{1/2}u_0}_2 + \norm{\phi^{1/2}u_{0x}}_2 + \norm{\phi^{1/2}u_{0xx}(t)}_2 + \norm{\phi^{1/2}\rho_0}_2 + \norm{\phi^{1/2}\rho_{0x}}_2\right)e^{CMt}.\nonumber
\end{align}
Let $\phi_n$ be as in the proof of Theorem~\ref{thm1}. Then \eqref{est.5} holds equally with $\phi$ replaced by $\phi_n$. By the definition of $F$ and \eqref{est.5}, there is a constant $\tilde C_1$ depending only on $\phi$ and $z_0$ such that $\norm{\phi_n F(u,\rho)}_1\leq\tilde C_1e^{2CMt}$. Using \eqref{Gx}, \eqref{young} and $ve^{-|\cdot|}\in L_p(\R)$ we conclude that $\norm{\phi_n(\partial_xG*F(u,\rho))}_p\leq\tilde C_2e^{2CMt}$ and, by \eqref{Gxx}, \eqref{young} and \eqref{phinF}, that
\begin{align}
\norm{\phi_n(\partial_x^2G*F(u,\rho))}_p&\leq\norm{\phi_n(G*F(u,\rho))}_p+\norm{\phi_nF(u,\rho)}_p\nonumber\\
&\leq\tilde C_3e^{2CMt}+M\left(\norm{\phi_nu}_p+\norm{\phi_nu_x}_p+\norm{\phi_n\rho}_p\right).\nonumber
\end{align}
Using \eqref{aux1}, \eqref{aux2} and \eqref{aux3} this yields
\begin{align}
& \frac{d}{dt}\left(\norm{\phi_nu}_p+\norm{\phi_nu_x}_p+\norm{\phi_nu_{xx}}_p+\norm{\phi_n\rho}_p+\norm{\phi_n\rho_x}_p\right) \nonumber\\
& \qquad\leq \tilde C_4M\left(\norm{\phi_nu}_p+\norm{\phi_nu_x}_p+\norm{\phi_nu_{xx}}_p+\norm{\phi_n\rho}_p+\norm{\phi_n\rho_x}_p\right)+\tilde C_5e^{2CMt} \nonumber
\end{align}
and the constants $\tilde C_j>0$, $j=1,\ldots,5$, depend only on $\phi$ and $z_0$. Integrating this equation and letting $n\to\infty$, we obtain the result of the corollary for $p\in[2,\infty)$. The case $p=\infty$ is again obtained from a standard limit argument. This achieves the proof.\hfill$\square$
\subsection{Proof of Corollary~\ref{cor2}} Conservation of \eqref{cond} follows from Theorem~\ref{thm1} with $p=\infty$ and $\phi(x)=\psi(x)=e^{|x|/2}(1+|x|)^{1/2}\left(\log(e+|x|)\right)^d$. Integration of \eqref{2CH} yields
$$
z(x,t)=z_0 - \begin{pmatrix}\int_0^t\partial_xG*F(u,\rho)(x,s)\d s+\int_0^t(uu_x)(x,s)\,ds \\ \int_0^t(\rho u_x+\rho_xu)(x,s)\,ds\end{pmatrix}.
$$
We now use that
\begin{align}
&\left|\int_0^t(uu_x)(x,s)\,ds\right|,\left|\int_0^t(\rho u_x)(x,s)\,ds\right|,\left|\int_0^t(\rho_xu)(x,s)\,ds\right|\nonumber\\
&\qquad\leq Cte^{-|x|}(1+|x|)^{-1}\left(\log(e+|x|)\right)^{-2d}\nonumber
\end{align}
and let
$$\Phi^\pm(t)=\frac{1}{2}\int_{-\infty}^\infty e^{\pm y}h(y,t)\d y,\quad h(x,t)=\frac{1}{t}\int_0^tF(u,\rho)(x,s)\d s.$$
By the arguments in \cite{Br12}, this achieves the asymptotic representation of $z$.\hfill$\square$

\end{document}